\documentclass[11pt]{amsproc}%
\usepackage{amssymb}
\usepackage{amsmath}
\usepackage{amsfonts}
\usepackage{graphicx}%
\providecommand{\U}[1]{\protect\rule{.1in}{.1in}}
\theoremstyle{plain}

\newtheorem{claim}{Claim}

\newtheorem{corollary}{Corollary}

\newtheorem{lemma}{Lemma}

\newtheorem{remark}{Remark}

\newtheorem{theorem}{Theorem}
\numberwithin{equation}{section}
\begin{document}
\title[Smooth extensions]{Smooth extensions of functions on separable Banach spaces}
\author{}
\author{D. Azagra}
\address{ICMAT (CSIC-UAM-UC3-UCM), Departamento de An{\'a}lisis Matem{\'a}tico,
Facultad Ciencias Matem{\'a}ticas, Universidad Complutense, 28040, Madrid, Spain}
\email{daniel\_azagra@mat.ucm.es}
\author{R. Fry}
\address{Department of Mathematics and Statistics, Thompson Rivers University,
Kamloops, B.C., Canada}
\email{rfry@stfx.ca }
\author{L. Keener}
\address{Department of Mathematics and Statistics, University of Northern British
Columbia, Prince Georege, B.C., Canada}
\email{keener@unbc.ca}
\keywords{Smooth extension, Banach space. The second named author is partly supported by
NSERC (Canada).}

\begin{abstract}
Let $X$ be a Banach space with a separable dual $X^{*}$. Let $Y\subset X$ be a
closed subspace, and $f:Y\rightarrow\mathbb{R}$ a $C^{1}$-smooth function.
Then we show there is a $C^{1}$ extension of $f$ to $X$.

\end{abstract}
\maketitle

\section{Introduction}

In this note we address the problem of the extension of smooth functions from
subsets of Banach spaces to smooth functions on the whole space. For our
results, smoothness is meant in the Fr{\'e}chet sense, and we shall restrict
our attention to real-valued functions. To state the problem more precisely,
given a Banach space $X$, a closed subset $Y,$ and a $C^{p}$-smooth function
$f:Y\rightarrow\mathbb{R},$ when is it possible to find a $C^{p}$-smooth map
$F:X\rightarrow\mathbb{R}$ such that $F\mid_{Y}=f$?

We should note that when $Y$ is a complemented subspace of an arbitrary Banach
space $X$, the extension problem can be easily solved. Indeed, let
$P:X\rightarrow Y$ be a continuous linear projection, and $f:Y\rightarrow
\mathbb{R}$ a $C^{p}$-smooth function. Then $F\left(  x\right)  =f\left(
Px\right)  $ defines a $C^{1}$ extension of $f$ to $X$. Unfortunately, not
every closed subspace $Y$ of a separable Banach space $X$ is complemented. In
fact, a classic result of Lindenstrauss and Tzafriri \cite{LT} states that the
only Banach space all of whose closed subspaces are complemented is (up to
renorming) a Hilbert space, so this trick only works when $X$ is a Hilbert space.

When $p=0,$ this question is the problem of the continuous extension of
functions from closed subsets. A complete characterization was given by the
well known theorem of Tietze (see e.g., \cite{Wi}) which we recall states that
$X$ is a normal space iff for every closed subset $Y\subset X$ and continuous
function $f:Y\rightarrow\mathbb{R},$ there exists a continuous extension
$F:X\rightarrow\mathbb{R}$ of $f.$

Such characterizations in the differentiable case, where $p\geq1,$ are more
delicate. When $X=\mathbb{R},$ $Y\subset X$ is a subset, and $p\geq1,$
necessary and sufficient conditions (in terms of divided differences) for the
existence of $C^{p}$-extensions to $\mathbb{R}$ of $C^{p}$ functions on $Y$
were given by H. Whitney \cite{W1}, \cite{W2}. Apparently, Whitney intended to
find such a characterization in the case $X=\mathbb{R}^{n}$ with $n>1,$ but a
sequel to the paper \cite{W2} never appeared.

\medskip

Major advances in this area occurred some twenty years later with the
fundamental work of Glaeser \cite{G} who solved the problem when $n\geq1$ and
$p=1$. Subsequent work included that of Brudnyi and Shvartsman \cite{BS1},
\cite{BS2}, Bierstone, Milman and Pawlucka \cite{BMP1}, \cite{BMP2}, and in
particular the striking results of C. Fefferman \cite{Fe1}, \cite{Fe2},
\cite{Fe3}. For example, in \cite{Fe2} a complete characterization is given of
when a real-valued function defined on a compact subset of $\mathbb{R}^{n}$ is
the restriction of a $C^{m}$-smooth map on $\mathbb{R}^{n}.$

\medskip

In this paper we consider the case when $X$ is a separable Banach space which
admits a $C^{1}$-smooth norm, a condition which is well known to be equivalent
to $X^{*}$ being separable \cite{DGZ}. Then if $Y\subset X$ is a closed
subspace and $f:Y\rightarrow\mathbb{R}$ is $C^{1}$-smooth, we show there
exists a $C^{1}$ extension $F:X\rightarrow\mathbb{R}$. If we require only that
$Y\subset X$ be closed and not necessarily a subspace, then a similar
conclusion holds under the stronger assumption that $f$ is defined on a
neighbourhood $U\supset Y$ and is $C^{1}$-smooth on $Y$ as a function on $X$
(i.e., $f^{\prime}\left(  y\right)  \in X^{\ast}$ for $y\in Y$ and
$y\rightarrow f^{\prime}\left(  y\right)  $ is continuous). We observe,
however, that in general the smooth extension problem has a negative solution.
We give here three examples.

\medskip

\begin{enumerate}
\item In \cite{Z} (see also \cite[Theorem II.8.3, page 82]{DGZ}) an example is
given of a separable Banach space $Y\subset X=C\left[  0,1\right]  ,$ and a
G{{\^a}}teaux smooth norm on $Y$ that cannot be extended to a G{{\^a}}teaux
smooth norm on $X.$

\medskip

\item Also in \cite{Z}, it is shown that for $1<p<2,$ there is a subspace
$Y\subset L_{p}$ isomorphic to Hilbert space, such that the Hilbertian norm of
$Y$ cannot be extended to a function $\varphi$ on $L_{p}$ which is Fr{\'e}chet
smooth on the unit sphere $S_{Y}$ of $Y$ as a function on $X$ with
$y\rightarrow\varphi^{\prime}\left(  y\right)  $ locally Lipschitz from
$S_{Y}$ to $X^{\ast}.$ Since every $C^{2}$ smooth function has a locally
Lipschitz derivative, this immediately shows that, given a $C^{\infty}$ smooth
function $f$ on $Y\subset L_{p}$ (with $1<p<2)$, in general there is no
$C^{2}$ smooth extension of $f$ to $L_{p}$.

\medskip

\item As suggested to us by R. Aron \cite{A} (see also Example 2.1 \cite{AB}).
Let $X=C\left[  0,1\right]  $ which has the Dunford-Pettis Property, and hence
the polynomial Dunford-Pettis Property (i.e., if $P:X\rightarrow\mathbb{R}$ is
a polynomial and $x_{i}\overset{w}{\rightarrow}0,$ then $P\left(
x_{i}\right)  \rightarrow P\left(  0\right)  $). Now $Y=l_{2}\subset X$ by the
Banach-Mazur Theorem, and we consider $f\left(  x\right)  =\left\Vert
x\right\Vert _{l_{2}}^{2}.$ If $f$ extended to a $C^{2}$-smooth function $F$
on an open neighbourhood of $l_{2}$ in $X,$ then
\[
P\left(  h\right)  =\left(  1/2\right)  F^{^{\prime\prime}}\left(  0\right)
\left(  h,h\right)
\]
would be a polynomial on $X.$ But $e_{i}\overset{w}{\rightarrow}0$ in $l_{2}$
and so in $X,$ and thus as noted above we would have $1=P\left(  e_{j}\right)
\rightarrow P\left(  0\right)  =0,$ a contradiction.
\end{enumerate}

\smallskip

One can compare the results of this note with the work of C.J. Atkin
\cite{At}. The programme of Atkin is to find smooth extension results for
smooth functions $f$ defined on finite unions of open, convex sets in
separable Banach spaces $X$ that do not admit $C^{p}$-smooth norms, or even
$C^{p}$-smooth bump functions. In order to achieve this, however, it is
assumed in \cite{At} that the function $f$ already possesses smooth extensions
to all of $X$ in a neighbourhood of every point in its domain. Finally we
mention the result \cite[Proposition VIII.3.8]{DGZ}, which states, in
particular, that for weakly compactly generated $X$ which admit $C^{p}$-smooth
bump functions, for any closed subset $Y\subset X$ and continuous function
$f:Y\rightarrow\mathbb{R},$ there exists a continuous extension of $f$ to $X$
which is $C^{p}$-smooth on $X\backslash Y.$

\smallskip

We remark that the situation for analytic maps is quite different. Indeed, the
paper by R. Aron and P. Berner \cite{AB} characterizes the existence of
analytic extensions from subspaces in terms of the existence of a linear
extension operator. In particular, in the real case they prove, among many
other equivalences, that if $Y$ is a closed subspace of a Banach space $X,$
the the following are equivalent (recall that $Z$ is a $\mathcal{C}$-space if
it is complemented in its second dual $Z^{\ast\ast})$:

\begin{enumerate}
\item For any $\mathcal{C}$-space $Z$, and (real) analytic map $f:Y\rightarrow
Z$ which is bounded on bounded sets, there exists an analytic extension
$F:X\rightarrow Z$ of $f,$ also bounded on bounded sets.

\item There exists a continuous, linear extension operator $T:Y^{\ast
}\rightarrow X^{\ast}.$
\end{enumerate}

\medskip

We have combined some recent work on smooth approximation of Lipschitz
mappings \cite{F} with some techniques of Moulis \cite{M}, the Bartle-Graves
selector theorem, and the classical method of Tietze to deduce our principal results.

\medskip

Our notation is standard, with $X$ typically denoting a (real) Banach space.
We shall denote an open ball with centre $x\in X$ and radius $r>0$ either by
$B_{r}\left(  x\right)  ,$ $B\left(  x;r\right)  ,$ or $B_{r}$ if the centre
is understood. We write the closed unit ball of a Banach space $X$ as $B_{X}.$
If $Y\subset X,$ we denote the restriction of a function $f:X\rightarrow
\mathbb{R}$ to $Y$ by $f\mid_{Y},$ and we say that a map $K:X\rightarrow
\mathbb{R}$ is an extension of $f:Y\rightarrow\mathbb{R}$ if $K\left(
y\right)  =f\left(  y\right)  $ for all $y\in Y.$ We denote the Fr{{\'e}}chet
derivative of a function $g$ at $x$ in the direction $h$ by $g^{\prime}\left(
x\right)  \left(  h\right)  .$ As noted above, if $Y\subset X$ and $U\supset
Y$ is open, we say that $f:U\rightarrow\mathbb{R}$ is $C^{1}$-smooth on $Y$ as
a function on $X$ if $f^{\prime}\left(  y\right)  \in X^{\ast}$ and
$y\rightarrow f^{\prime}\left(  y\right)  $ is continuous for $y\in Y.$ For
any undefined terms we refer the reader to \cite{FHHMPZ}, \cite{DGZ}.

\medskip

\section{Main Results}

\noindent An essential tool shall be the following consequence of
the main theorem in \cite{F}, see also \cite{HJ}.
\begin{lemma}
\label{Main} There exists a constant $C_{0}\geq1$ such that, for
every separable Banach $X$ space with a $C^{1}$-smooth norm, for
every subspace $Y\subseteq X$, every Lipschitz function
$f:X\rightarrow\mathbb{R}$, and every $\varepsilon>0$, there
exists a $C^{1} $-smooth function $K:X\rightarrow\mathbb{R}$ such
that
\begin{enumerate}
\item $\left\vert f(x)-K(x)\right\vert <\varepsilon$ for all $x\in
X$,
\item $\textrm{Lip}(K)\leq C_0 \textrm{Lip}(f)$, and
\item $\|K'(y)\|_{X^{*}}\leq C_0 \textrm{Lip}(f_{|_Y})$ for all $y\in Y$ (in
particular the Lipschitz constant of the restriction of $K$ to $Y$
is of the order of the Lipschitz constant of the restriction of
$f$ to $Y$).
\end{enumerate}
\end{lemma}
\noindent This lemma can be deduced with some work from the
results in either \cite{F,HJ} but here, for the sake of
completeness, we shall give a self-contained proof which moreover
provides a simple method of constructing {\em sup-partitions of
unity}.

\medskip

\noindent {\bf Proof of the Lemma.} Let us first assume that
$f:X\to [1,1001]$. Define $\eta=\textrm{Lip}(f_{|_Y})$,
$L=\textrm{Lip}(f)$, $R=1/\eta$ and $r=1/L$ (in the event that
$\textrm{Lip}(f_{|_Y})=0$, take any $\eta\in (0,L)$, and observe
that when $\textrm{Lip}(f)=0$ the result is trivial, so we may
assume $L>0$). Obviously, $\eta\leq L$ and $r\leq R$.

Since $f$ is $\eta$-Lipschitz on $Y$ and $Y$ is separable we can
cover $Y$ by a countable family of balls $B(y_n, R)$ of radius
$R$, where $\{y_n\}$ is a dense subset in $Y$, in such a way that
if $y, y'\in B(y_n, 4R)\cap Y$ then $|f(y)-f(y')|\leq 8$.
Similarly, since $f$ is $L$-Lipschitz on $X$, we can cover the set
$\{x\in X: \textrm{dist}(x, Y)\geq r/4\}$ by a countable family of
balls $B(x_n, r/32)$ of radius $r/32$, where $\{x_n\}$ is dense in
$\{x\in X: \textrm{dist}(x, Y)\geq r/4\}$, with the properties
that the balls $B(x_n, r/8)$ of radius $4r/32$ do not touch the
set $\{x\in X: \textrm{dist}(x, Y)< r/8\}$, and that if $x, x'\in
B(x_n, r/8)$ then $|f(x)-f(x')|\leq 1/4$.

Also note that the open slabs $D_{y_n}:=\{x\in X: \textrm{dist}(x,
Y)< r, \|x-y_n\|< R\}$ cover $Y$ and if we denote
$$D_{y_n}^{4}:=\{x\in X: \textrm{dist}(x, Y)< r, \|x-y_n\|< 4R\}$$
then these sets have the property that if $x, x'\in D_{y_n}^4$
then $|f(x)-f(x')|\leq 10$.

\begin{claim}
\label{construction of the sup partition of unity} There exists a
sequence of $C^1$ functions $\varphi_{n}
:X\to\mathbb{R}$ with the following properties:
\begin{enumerate}
\item The collection $\{\varphi_{n} :X\to[0,1] \, | \, n\in\mathbb{N}\}$ is
uniformly Lipschitz on $X$, with Lipschitz constant $8/r=8L$.
\item $\|{\varphi_n}'(y)\|_{X^{*}}\leq 2/R=2\eta$ for all $y\in Y$. In fact,
$$
\textrm{Lip}({\varphi_{n}}_{|_{\{x\in X \, : \, \text{ \em
dist}(x, Y)<r/4\}}})\leq 2\eta.$$
\item For each $x$ with $\textrm{dist}(x,Y)<r/4$ there exists $n\in\mathbb{N}$ with
$\varphi_{n}(x)=1$.
\item For each $x\in X$ there exists $\delta>0$ and $n_{x}\in\mathbb{N}$ such
that for $z\in B(x, \delta)$ and $n>n_{x}$ we have $\varphi_{n}
(z)=0$.
\item $\varphi_{n}(x)=0$ for all $x\notin D_{y_n}^{4}$.
\end{enumerate}
\end{claim}
The existence of a family of functions $\varphi_n$ satisfying
properties $(1), \, (3), \, (4)$ and $(5)$ on all of $X$ is known
from \cite{F, HJ}. What is new about this claim is that, in the
present situation, one can also require (property $(2)$) that the
derivatives of $\varphi_n$ are bounded on $Y$ by a constant of the
order of $\textrm{Lip}(f_{|_Y})$, which could be very small
compared to the global Lipschitz constant of $f$. We say that the
collection of functions $\varphi_n$ forms a {\em sup-partition of
unity} on $\{x\in X: \textrm{dist}(x, Y)<r\}$, subordinated to the
covering $\{D_{y_n}: n\in\mathbb{N}\}$.

\medskip

\noindent {\bf Proof of the Claim.} Define subsets
$A_{1}=\{u_{1}\in\mathbb{R}: -1\leq u_{1}\leq 4r\}$, and, for
$n\geq2$,
\[
A_{n} =\{\{u_{j}\}_{j=1}^{n}\in\ell^{n}_{\infty} : -1-R\leq
u_{n}\leq4R, \, 2R\leq u_{j}\leq M_{n}+2 \text{ for } 1\leq j\leq
n-1\},
\]
\[
A^{\prime}_{n}=\{\{u_{j}\}_{j=1}^{n}\in\ell^{n}_{\infty} : -1\leq u_{n}%
\leq3R, \, 3R\leq u_{j}\leq M_{n}+2-R \text{ for } 1\leq j\leq
n-1\},
\]
\[
\text{where } M_{n}=\sup\left\{ \| x-y_{j}\| :x\in B(y_{n}, 4R),\
1\leq j\leq n\right\} .
\]

Let $b_{n}:\ell_{\infty}^{n}\rightarrow [0,2]$ be the function
defined by
\[
b_{n}(y)=\max\{0,1-\frac{1}{R}\text{dist}_{\infty}(y,A_{n}^{\prime})\},
\]
where $\text{dist}_{\infty}(y,A)=\inf\{\Vert
y-a\Vert_{\infty}:a\in A\}$. It is clear that
support$(b_{n})=A_{n}$, that $b_{n}=1$ on $A_{n}^{\prime}$, and
that $b_{n}$ is $(1/R)$-Lipschitz (note in particular that the
Lipschitz constant of $b_{n}$ does not depend on $n$).

Since the function $b_{n}$ is uniformly continuous and bounded on
$\mathbb{R}^{n}$, it is a standard fact that the normalized
integral convolutions of $b_{n}$ with the Gaussian-like kernels
$y\mapsto G_{\kappa }(y):=
e^{-\kappa\sum_{j=1}^{{n}}2^{-j}y_{j}^{2}}$,
\[
x\mapsto\frac{1}{T_{\kappa}}b_{n} * G_{\kappa}(x)=\frac{1}{\int_{\mathbb{R}%
^{n}}e^{-\kappa\sum_{j=1}^{{n}}2^{-j}y_{j}^{2}}dy} \int_{\mathbb{R}^{n}}%
b_{n}(y)e^{-\kappa\sum_{j=1}^{n}2^{-j}(x_{j}-y_{j})^{2}}dy,
\]
\[
\text{where } T_{\kappa}=\int_{\mathbb{R}^{n}}e^{-\kappa\sum_{j=1}^{{n}}%
2^{-j}y_{j}^{2}}dy,
\]
converge to $b_{n}$ uniformly on $\mathbb{R}^{n}$ as
$\kappa\to+\infty$. Therefore, for each $n\in\mathbb{N}$ we can
find $\kappa_{n}>0$ large enough so that
\[
|b_{n}(x)-\frac{1}{T_{\kappa_{n}}}b_{n} * G_{\kappa_{n}}(x)|\leq
1/10 \,\, \text{ for all } \,\, x\in\mathbb{R}^{n}. \eqno(*)
\]

Define $\nu_{n}:\ell^{n}_{\infty}\to\mathbb{R}$ by
\[
\nu_{n} (x) := \frac{1}{T_{\kappa_{n}}}b_{n} * G_{\kappa_{n}}(x) =
\frac
{1}{T_{n}}\int_{\mathbb{R}^{n}}b_{n}(y)e^{-\kappa_{n}\sum_{j=1}^{n}%
2^{-j}(x_{j}-y_{j})^{2}}dy.
\]
Let us note that
\[
\frac{1}{T_{n}}\int_{\mathbb{R}^{n}}b_{n}(y)e^{-\kappa_{n}%
\sum_{j=1}^{n}2^{-j}(x_{j}-y_{j})^{2}}dy= \frac{1}{T_{n}}\int_{\mathbb{R}^{n}%
}b_{n}(x-y)e^{-\kappa_{n}\sum_{j=1}^{n}2^{-j}y_{j}^{2}}dy,
\]
and so
\begin{align*}
\left\vert \nu_{n}\left(  x\right)  -\nu_{n}\left(
x^{\prime}\right) \right\vert  &  =\left\vert
\frac{1}{T_{n}}\int_{\mathbb{R}^{n}}\left( b_{n}\left(  x-y\right)
-b_{n}\left(  x^{\prime}-y\right)  \right)
e^{-\kappa_{n}\sum_{j=1}^{n}2^{-j}y_{j}^{2}}dy\right\vert \\
& \\
&  \leq\frac{1}{T_{n}}\int_{\mathbb{R}^{n}}\left\vert b_{n}\left(
x-y\right)
-b_{n}\left(  x^{\prime}-y\right)  \right\vert e^{-\kappa_{n}\sum_{j=1}%
^{n}2^{-j}y_{j}^{2}}dy\\
& \\
&  \leq\frac{1}{R}\ \left\Vert x-x^{\prime}\right\Vert _{\infty}\frac{1}%
{T_{n}}\int_{\mathbb{R}^{n}}e^{-\kappa_{n}\sum_{j=1}^{n}2^{-j}y_{j}^{2}}dy\\
& \\
&  =\frac{1}{R}\left\Vert x-x^{\prime}\right\Vert _{\infty}.
\end{align*}
Hence, $\nu_{n}$ is $\frac{1}{R}$-Lipschitz. Note also that $0\leq
\nu_n (x)\leq\|b_n\|_{\infty}=1$ for all $x\in X$.

\medskip

Now take a $C^{\infty}$ function $\alpha:\mathbb{R}\to [0,1]$ such
that:
\begin{enumerate}
\item[{}] $\alpha$ is $2$-Lipschitz;
\item[{}] $\alpha(t)=0$ if $t\leq 1/10$;
\item[{}] $\alpha(t)=1$ if $t\geq 9/10$.
\end{enumerate}
Then the composition $\alpha\circ\nu_n$ is a $C^{\infty}$ function
so that
\begin{enumerate}
\item[{}] $\alpha\circ\nu_n$ is $2/R$-Lipschitz;
\item[{}] $\alpha\circ\nu_n(x)=0$ if $x\notin A_n$;
\item[{}] $\alpha\circ\nu_n=1$ if $x\in {A_n}'$.
\end{enumerate}

Consider the quotient space $X/Y$, with its quotient map $q:X\to
X/Y$. The mapping $T\ni (X/Y)^{*}\to T\circ q \in X^{*}$ defines a
continuous linear injection from $(X/Y)^{*}$ into $X^{*}$, and
since $X^{*}$ is separable so is $(X/Y)^{*}$. Hence $X/Y$ has an
equivalent $C^1$ smooth norm (which we will also denote
$\|\cdot\|$) with the property that
    $$
\textrm{dist}(x, Y)\leq \|q(x)\|\leq 2\textrm{dist}(x,Y) \textrm{
for all } x\in X.
    $$
In particular the function $x\mapsto \|q(x)\|$ is $2$-Lipschitz on
$X$, as is easily checked.

Take also a $C^\infty$ function $\beta:\mathbb{R}\to [0,1]$ such
that
\begin{enumerate}
\item[{}] $\beta$ is $3/r$-Lipschitz;
\item[{}] $\beta(t)=0$ if $t\geq r$;
\item[{}] $\beta(t)=1$ if $t\leq r/2$.
\end{enumerate}

\noindent Next, consider the map $\lambda_{n}:X\rightarrow
l_{\infty}^{n}$ given by
\[
\lambda_{n}\left(  x\right)  =\left( \| x- y_{1}\| ,..., \| x-
y_{n}\| \right).
\]

\medskip

\noindent Then for $n\geq1$ we define the maps $\varphi
_{n}:X\rightarrow\mathbb{R}$ by
\begin{align*}
\varphi_{n}\left(  x\right)   &  = \, \beta\left(\|q(x)\|\right)\,
\alpha\left(\nu_{n}\left( \lambda_{n}\left( x\right)
\right)\right) \, = \, \beta\left(\|q(x)\|\right)\,
\alpha\left(\nu_{n}\left( \left\{ \| x-y_{j} \| \right\}
_{j=1}^{n}\right)\right).
\end{align*}
Since $\nu_n$ is constant in a neighborhood of each point $v=(v_1,
..., v_n)\in \mathbb{R}^{n}$ with $v_i=0$ for some $i$, it is
immediately seen that $\varphi_n$ is of class $C^1$ on $X$.

Now, if $\textrm{dist}(x, Y)< r/4 >\textrm{dist}(x', Y)$ we have
that $\beta(\|q(x)\|)=1=\beta(\|q(x')\|)$ and
\begin{eqnarray*}
& &\left\vert \varphi_{n}\left(  x\right)  -\varphi_{n}\left(
x^{\prime}\right) \right\vert  =\left\vert
\alpha\circ\nu_{n}\left( \lambda_{n}\left(  x\right) \right)
-\alpha\circ\nu_{n}\left( \lambda_{n}\left(  x^{\prime}\right)
\right) \right\vert \\
& &\leq\frac{2}{R}\left\Vert \lambda_{n}\left(  x\right)
-\lambda_{n}\left( x^{\prime}\right)  \right\Vert _{\infty} =
\frac{2}{R}\left\Vert \left\{  \|  x-x_{j} \| -\|
x^{\prime}-x_{j} \|  \right\}  _{j=1}^{n}\right\Vert _{\infty}\\
& &  \leq\frac{2}{R}\left\Vert x-x^{\prime}\right\Vert _{X},
\end{eqnarray*}
hence the collection $\left\{  \varphi_{n}\right\}  $ is uniformly
Lipschitz on the open neighborhood $\{x\in X: \textrm{dist}(x,
Y)<r/4\}$ of the subspace $Y$, with constant $\frac{2}{R}=2\eta$.
In particular we have that
    $$
\|{\varphi_n}'(y)\|_{X^{*}}\leq 2\eta \, \textrm{ for all } y\in
Y,
    $$
which shows $(2)$. On the other hand, from the definition of the
$\varphi_n$, it is immediately checked that these functions are
uniformly Lipschitz on all of $X$, with constant $2/R+6/r\leq
8/r$. This shows $(1)$.

\medskip

Let us show (3). For each fixed $x\in X$ with $\textrm{dist}(x,
Y)<r/4$ there exists $n_{x}$ with $x\in B(y_{n_x},3R)$ but with
$x\notin B(y_i,3R)$ for $i<n_{x}$. This implies that the point
$\left( \| x-y_{1}\|  , \| x-y_{2} \| ,..., \| x-y_{n_{x}} \|
\right) $ belongs to $A_{n_{x}}^{\prime}$, where the function
$\alpha\circ\nu_{n_{x}}$ takes the value $1$. Besides
$\beta(\|q(x)\|)=1$. Hence by the definition of $\varphi_n$, we
have $\varphi_{n_{x}}(x)=1$.

\medskip

Property (5) is shown similarly: if $\|x-y_{n}\|\geq 4R$ then the
point $(\|x-y_{1}\|, ..., \|x-y_{n}\|)$ lies in a region of
$\mathbb{R}^{n}$ where the function $\alpha\circ\nu_{n}$ takes the
value $0$, hence $\varphi_n (x)=0$. Or, if $\textrm{dist}(x,
Y)\geq r$ then $\|q(x)\|\geq r$ and $\beta(\|q(x)\|)=0$, hence
$\varphi_n (x)=0$.

\medskip

We finally show (4). If $\textrm{dist}(x,Y)\leq r\leq R$ then,
since the sequence $\{y_n\}$ is dense in $Y$, there exists
$n_x\in\mathbb{N}$ such that $\|x-y_{n_x}\|<2R$. Take
$\delta=2R-\|x-y_{n_x}\|>0$. Then for all $z\in B(x, \delta)$ we
also have $\|z-y_{n_x}\|<2R$ and, by the definition of $A_n$,
    $$
\lambda_n (z)=\left( \|z-y_1\|, ..., \|z-y_{n_x}\|, ...,
\|z-y_{n}\|\right)\notin A_n \, \textrm{ for } n>n_x,
    $$
hence, bearing in mind that $\alpha\circ\nu_x = 0$ outside $A_n$,
we get
    $$
\varphi_n(z)=0 \, \textrm{for all } n>n_x, z\in B(x, \delta).
    $$
On the other hand, if $\textrm{dist}(x,Y)>r$ then
$\beta(\|q(z)\|)=0$ for all $z\in B(x, \delta')$, where
$\delta'=\textrm{dist}(x, Y)-r>0$, and therefore $\varphi_n(z)=0$
for all $n\in\mathbb{N}, z\in B(x, \delta')$. \,\, $\blacksquare$

\bigskip

\begin{remark}
{\em Note that if $Y=X$ then $q=0$ and $\beta(\|q(x)\|)=1$ for all
$x$ (hence there is no need to use this term in the definition of
$\varphi_n$). In this case the above proof gives a simple method
of constructing $2/R$-Lipschitz sup-partitions of unity
subordinated to any covering by balls of radius $R$ of $X$.
Moreover these sup-partitions of unity are of the same order of
smoothness as the norm of $X$.}
\end{remark}

\bigskip

By replacing $y_n$ with $x_n$, $R$ with $r/32$, and $\beta$ with a
different $C^1$ function $\beta:\mathbb{R}\to [0,1]$ such that
\begin{enumerate}
\item[{}] $\beta(t)=0$ for $t\leq r/4 - r/16$
\item[{}] $\beta(t)=1$ for $t\geq r/4 - r/32$
\item[{}] $\textrm{Lip}(\beta)\leq 36/r$,
\end{enumerate}
one can similarly show:
\begin{claim}
There exists a sequence of $C^1$ functions $\psi_{n}
:X\to\mathbb{R}$ with the following properties:
\begin{enumerate}
\item The collection $\{\psi_{n} :X\to[0,1] \, | \, n\in\mathbb{N}\}$ is
uniformly Lipschitz on $X$, with Lipschitz constant $136/r=136L$.
\item For each $x$ with $\textrm{dist}(x,Y)\geq r/4$ there exists $n\in\mathbb{N}$ with
$\psi_{n}(x)=1$.
\item For each $x\in X$ there exists $\delta>0$ and $n_{x}\in\mathbb{N}$ such
that for $z\in B(x, \delta)$ and $n>n_{x}$ we have $\varphi_{n}
(z)=0$.
\item $\psi_{n}(x)=0$ for all $x\notin B(x_n, r/8)$.
\end{enumerate}
In particular all of the functions $\psi_n$ vanish on the set
$\{x\in X : \textrm{dist}(x, Y)<r/8\}$.
\end{claim}

Now let $\|\cdot\|_{c_0}$ be a $C^{\infty}$ smooth equivalent norm
to the usual norm $\|\cdot\|_{\infty}$ of $c_0$ and such that
    $$
\|x\|_{\infty}\leq\|x\|_{c_0}\leq 2\|x\|_{\infty} \textrm{ for all
} x\in c_0.
    $$
Let us define a collection of $C^1$ functions $\Phi_n:X\to [0,1]$
by
    $$
\Phi_n(x)=
  \begin{cases}
    \varphi_{k}(x) & \text{ if } n=2k-1 \text{ is odd },  \\
    \psi_{k}(x) & \text{ if } n=2k \text{ is even }.
  \end{cases}
    $$
Notice that, according to properties $(4)$ of Claim 1 and $(3)$ of
Claim 2, the mapping $X\ni x\to \{\Phi_{n}(x)\}_{n=1}^{\infty}\in
c_0$ is well defined and $C^1$ smooth (as the tails of the
sequence eventually vanish locally).

Define a function $g:X\to\mathbb{R}$ by
    $$
g(x)=\frac{\|\{a_n\Phi_{n}(x)\}_{n=1}^{\infty}\|_{c_0}}{\|\{\Phi_{n}(x)\}_{n=1}^{\infty}\|_{c_0}},
    $$
where
    $$
a_n=
  \begin{cases}
    f(y_k) & \text{ if } n=2k-1 \text{ is odd },  \\
    f(x_k) & \text{ if } n=2k \text{ is even }.
  \end{cases}
    $$
The function $g$ is well defined because
$\|\{\Phi_{n}(x)\}_{n=1}^{\infty}\|_{c_0}\geq
\|\{\Phi_{n}(x)\}_{n=1}^{\infty}\|_{\infty}=1$ by properties $(3)$
of Claim 1 and $(2)$ of Claim 2, and is $C^1$ smooth on $X$ by the
previous observation and because $a_n\geq 1$.

Since the functions $\Phi_n$ are $136\times L$-Lipschitz and
$|a_n|\leq 1001$ we have
\begin{eqnarray*}
& &
\|\{a_n\Phi_{n}(x)\}_{n=1}^{\infty}-\{a_n\Phi_{n}(z)\}_{n=1}^{\infty}\|_{c_0}\leq
2\|\{a_n\left(\Phi_{n}(x)-\Phi_{n}(z)\right)\}_{n=1}^{\infty}\|_{\infty}\leq\\
& & 2002\times 136\times L\|x-y\|,
\end{eqnarray*}
that is the function $\|\{a_n\Phi_n
(\cdot)\}_{n=1}^{\infty}\|_{c_0}$ is $2002\times 136\times
L$-Lipschitz on $X$, and is bounded by $2002$. Similarly, since
the function $t\mapsto 1/t$ is $1$-Lipschitz on $[1, \infty)$ and
$\{\Phi_n (\cdot)\}_{n=1}^{\infty}$ is bounded below by $1$, we
have that the function $1/\|\{\Phi_n
(\cdot)\}_{n=1}^{\infty}\|_{c_0}$ is $1\times 2\times 136\times
L$-Lipschitz on $X$ and bounded above by $1$. Therefore the
product satisfies
    $$
\textrm{Lip}(g)\leq 2002\times (1\times 2\times 136\times L) \, +
\, 1\times (2002\times 136\times L)=816816\times L.
    $$
When we restrict $g$ to the set $\{x\in X : \textrm{dist}(x,
Y)<r/8\}$, all the even terms of the sequence
$\{\Phi_{n}(x)\}_{n=1}^{\infty}$ vanish, so the only functions
that matter are the $\varphi_k$, which are $2\eta$-Lipschitz on
this set, and the above calculation can be performed replacing $L$
with $\eta$ to show that
    $$
\textrm{Lip}(g_{|_{\{x\in X \, : \, \textrm{ dist}(x,
Y)<r/4\}}})\leq 816816 \times \eta,
    $$
which implies
    $$
\|g'(y)\|_{X^{*}}\leq 816816\eta \textrm{ for all } y\in Y.
    $$
Finally, bearing in mind that the supports of the $\varphi_n$ are
contained in the slabs $D_{y_n}^{4}$, that the supports of the
$\psi_n$ are contained in the balls $B(x_n, r/8)$, and that on
each of these sets the oscillation of $f$ is bounded by $10$, it
is easy to check that
    $$
|f(x)-g(x)|\leq 20 \textrm{ for all } x\in X.
    $$

This argument proves the Lemma in the case when $\varepsilon=20$
and $f:X\to [0, 1000]$.

We next see that this result remains true for functions $f$ taking
values in $\mathbb{R}$ if we replace $20$ with $50$ and we allow
$C_0$ to be slightly larger than $816816$. Indeed, by considering
the function $h=\theta\circ g$, where $\theta$ is a $C^{\infty}$
smooth function $\theta:\mathbb{R}\to[0,1000]$ such that
$|t-\theta(t)|\leq 30$ if $t\in[0, 1000]$, $\theta(t)=0$ for
$t\leq 21$, and $\theta(t)=1000$ for $t\geq 979$, we get the
following result: there exists $C_{0}:=816816\times
\text{Lip}(\theta)$ such that for every $L$-Lipschitz function
$f:X\to[0,1000]$ whose restriction to $Y$ is $\eta$-Lipschitz
there exists a $C^{1}$ function $h:X\to[0,1000]$ such that

\begin{enumerate}
\item $|f(x)-h(x)|\leq 50$ for all $x\in X$

\item $h$ is $C_{0} L$-Lipschitz

\item $\|h'(y)\|_{X^{*}}\leq C_0 \eta$

\item $f(x)=0 \implies h(x)=0$, and $f(y)=1000 \implies h(y)=1000$.
\end{enumerate}

Now, for a $L$-Lipschitz function
$f:X\rightarrow\lbrack0,+\infty)$ so that
$\textrm{Lip}(f_{|_Y})=\eta$, we can write
$g(x)=\sum_{n=0}^{\infty}f_{n}(x)$, where
\[
f_{n}(x)=%
\begin{cases}
f(x)-1000n & \text{ if }1000n\leq f(x)\leq1000(n+1),\\
0 & \text{ if }f(x)\leq1000n,\\
1000 & \text{ if }1000(n+1)\leq f(x)
\end{cases}
\]
and the sum is locally finite. The functions $g_{n}$ are clearly
$L$-Lipschitz, satisfy $\textrm{Lip}((g_n)_{|_Y})\leq\eta$ and
take values in the interval $[0,1000]$, so there are $C^{1}$
functions $h_{n}:X\rightarrow\lbrack0,1000]$ such that for all
$n\in\mathbb{N}$ we have that $h_{n}$ is $C_{0} L$-Lipschitz,
$\|h'_n(y)\|_{X^{*}}\leq C_0 \eta$ for all $y\in Y$,
$|f_{n}-h_{n}|\leq 50$, and $h_{n}$ is $0$ or $1000$ wherever
$f_{n}$ is $0$ or $1000$. It is easy to check that the function
$h:X\rightarrow\lbrack0,+\infty)$ defined by
$h=\sum_{n=0}^{\infty}h_{n}$ is $C^{1}$ smooth, $C_{0}$-Lipschitz,
and satisfies $|f-h|\leq 50$ and $\|h'(y)\|_{X^{*}}\leq C_0 \eta$.
This argument shows that there is $C_{0}\geq1$ such that for any
$L$-Lipschitz function $f:X\rightarrow\mathbb{[}0,+\infty)$ with
$\textrm{Lip}(f)=\eta$, there exists a $C^{1}$ function
$h:X\rightarrow\lbrack0,+\infty)$ such that

\begin{enumerate}
\item $|f(x)-h(x)|\leq 50$ for all $x\in X$

\item $h$ is $C_{0}$-Lipschitz

\item $\|h'(y)\|_{X^{*}}\leq C_0 \eta$ for all $y\in Y$

\item $f(x)=0 \implies h(x)=0$.
\end{enumerate}

Finally, for an arbitrary $L$-Lipschitz function
$f:X\to\mathbb{R}$, we can
write $f=f^{+} - f^{-}$ and apply this result to find $C^{1}$ smooth, $C_{0}%
$-Lipschitz functions $h^{+}, h^{-} :X\to[0, +\infty)$ so that
$h:=h^{+}- h^{-}$ is $C^{1}$ smooth, $C_{0} L$-Lipschitz,
$\|h'\|_{X^{*}}\leq C_0 \eta$ on $Y$, $|f-h|\leq 50$. This proves
the Lemma for $\varepsilon=50$.

For an arbitrary $\varepsilon \in (0,50)$, let us consider the
function $g:X\rightarrow\mathbb{R}$ defined by
$g(x)=\frac{50}{\varepsilon}f(\frac{\varepsilon}{50} x)$. It is
immediately checked that $\textrm{Lip}(g)=\textrm{Lip}(f)=L$ and
$\textrm{Lip}(g_{|_Y})=\textrm{Lip}(f_{|_Y})=\eta$, so by the
result above there exists a $C^{1}$ smooth, $C_{0} L$-Lipschitz
function $h$ with $\|h'\|_{X^{*}}$ bounded by $C_0 \eta$ on $Y$
and such that $|g(x)-h(x)|\leq 50$ for all $x$, which implies that
the function $K(z):=\frac{\varepsilon}{50} h(\frac{50
}{\varepsilon}z)$ is $C_{0}\eta$-Lipschitz and satisfies
$|f(z)-K(z)|\leq \varepsilon$ for all $z\in X$. $\blacksquare$

\bigskip

\noindent We next establish the existence of a continuous and bounded
selection of the Hahn-Banach extension operator $y^{\ast}\in Y^{\ast
}\rightarrow G(y^{\ast})\in2^{X^{\ast}},$ where%
\[
G(y^{\ast})=\{x^{\ast}\in X^{\ast}:x^{\ast}(y)=y^{\ast}(y)\text{ for all }y\in
Y\}.
\]

\begin{lemma}
\label{existence of extension operators} For every Banach space $X$ and every
closed subspace $Y\subset X$ there exist a continuous mapping $H:Y^{\ast
}\rightarrow X^{\ast}$ and a number $M\geq1$ such that

\begin{enumerate}
\item $H(y^{*})(y)=y^{*}(y)$ for every $y^{*}\in Y^{*}$, $y\in Y$;

\item $\|H(y^{*})\|_{X^{*}}\leq M\|y^{*}\|_{Y^{*}}$ for every $y^{*}\in Y^{*}
$.
\end{enumerate}
\end{lemma}

\medskip

\noindent\textbf{Proof.\ \ }This is a consequence of the Bartle-Graves
selector theorem (see \cite[page 299]{DGZ}) which states: \emph{Let }$W$\emph{
and }$Z$\emph{ be Banach spaces and let }$T$\emph{ be a bounded linear mapping
of }$W$\emph{ onto }$Z$\emph{. Then there exists a continuous (nonlinear in
general) mapping }$B$\emph{ of }$Z$\emph{ into }$W$\emph{ such that}\textit{
}$\left(  T\circ B\right)  w=w$\textit{ for every }$w\in W$\textit{.}
Moreover, it follows from the proof of this result that there exists an $M>1$
such that $\left\Vert B(w)\right\Vert \leq M\Vert w\Vert$. If we apply this
theorem with $W=X^{\ast}$, $Z=Y^{\ast}$, to the mapping $T:X^{\ast}\rightarrow
Y^{\ast}$ defined by $T(x^{\ast})=x_{|_{Y}}^{\ast}$, which is a continuous
linear surjection with $\Vert T\Vert=1$ (by the Hahn-Banach theorem), we
obtain our continuous map $H=B:Y^{\ast}\rightarrow X^{\ast}$ with the property
that the the restriction of $H(y^{\ast})$ to $Y$ is $y^{\ast}$, for every
$y^{\ast}\in Y^{\ast}$, and such that $\Vert H(y^{\ast})\Vert_{X^{\ast}}\leq
M\Vert y^{\ast}\Vert_{Y^{\ast}}$. $\blacksquare$

\medskip

Now we are in a situation to deduce an approximation result which is of
independent interest and which, combined with some ideas of the Tietze proof,
will yield our main results on smooth extension.

\begin{theorem}
Let $X$ be a separable Banach space which admits a $C^{1}$-smooth norm, and
$Y\subset X$ a closed subspace. Let $f:Y\rightarrow\mathbb{R}$ be a $C^{1}%
$-smooth function, and $F$ a continuous extension of $f$ to $X$. Let
$H:Y^{\ast}\rightarrow X^{\ast}$ be any extension operator as in Lemma
\ref{existence of extension operators}. Then, for every $\varepsilon>0$, there
exists a $C^{1}$-smooth map $g:X\rightarrow\mathbb{R}$ such that

\begin{enumerate}
\item $\left\vert F\left(  x\right)  -g\left(  x\right)  \right\vert
<\varepsilon$ on $X,$ and

\item $\left\Vert H(f^{\prime}\left(  y\right)  )-g^{\prime}\left(  y\right)
\right\Vert _{X^{\ast}}<\varepsilon$ on $Y$.
\end{enumerate}

Furthermore, if the given $C^{1}$ function $f$ is Lipschitz on $Y$ and $F$ is
a Lipschitz extension of $f$ to $X$ with $\text{Lip}(F)=\text{Lip}(f)$ (for
instance $F(x)=\inf_{y\in Y}\{f(y)+\text{Lip}(f)\Vert x-y\Vert\}$), then the
function $g$ can be chosen to be Lipschitz on $X$ and with the additional
property that

\begin{enumerate}
\item[(3)] $\text{Lip}(g)\leq C \text{Lip}(f),$
\end{enumerate}

where $C >1$ is a constant only depending on $X$.
\end{theorem}

\noindent\textbf{Proof.\ \ }First note that by the Tietze Theorem, the
continuous extension $F$ always exists. We modify the proof of Theorem 4 in
\cite{AFGJL} employing Lemma 1. It will be convenient to use the following
notation: given a point $y_{k}\in Y,$ we define $T_{k}$ to be the natural
$H$-extension of the first order Taylor Polynomial of $f$ at $y_{k}$; namely,
$T_{k}\left(  x\right)  =f\left(  y_{k}\right)  +H(f^{\prime}\left(
y_{k}\right)  )\left(  x-y_{k}\right)  $. Note in particular that $T_{k}\in
C^{\infty}(X,\mathbb{R})$, with $T_{k}^{\prime}(x)=H(f^{\prime}(y_{k}))$ for
every $x\in X$, and $T_{k}^{\prime}(y)\mid_{Y}=f^{\prime}(y_{k})$ for all
$y\in Y$.

\medskip

\noindent Now, using the separability of $X,$ the closedness of $Y\subset X,$
and the continuity of $F,$ we can construct a covering $\mathcal{C=}\left\{
B_{r_{j}}\right\}  _{j=1}^{\infty}\cup\left\{  B_{s_{k}}\right\}
_{k=1}^{\infty}$ of $X$, by open balls with centres $x_{j}$ and $y_{k}$
respectively, with the following properties:

\medskip

(i). We have $B_{2r_{j}}\subset X\backslash Y,$ and $|F\left(  x\right)
-F\left(  x_{j}\right)  |<\varepsilon/2C_{0}$ on $B_{2r_{j}}$,

\medskip

(ii). The collection $\left\{  B_{s_{k}}\right\}  _{k}\subset X$ covers $Y$
with centres $y_{k}\in Y$ and radii $s_{k}$ chosen using the smoothness of $f$
on $Y$ and the norm-norm continuity of the extension operator $H$, so that
$\left\Vert T_{k}^{\prime}\left(  y\right)  -f^{\prime}\left(  y\right)
\right\Vert _{Y^{\ast}}<\varepsilon/8C_{0}$ and $\left\Vert T_{k}^{\prime
}\left(  y\right)  -H(f^{\prime}(y))\right\Vert _{X^{\ast}}<\varepsilon
/8C_{0}$ on $B_{2s_{k}}\cap Y.$


\medskip

\noindent It will be useful in the sequel, to employ an alternate notation for
the open balls $B_{r_{j}}$ and $B_{s_{k}}.$ We let $\beta:\mathbb{N}%
\rightarrow\mathcal{C}$ be a bijection where for each $i$, $\beta
(i)=B(\beta_{1}(i);\beta_{2}(i))$. Let $\varphi_{j}\in C^{1}\left(  X,\left[
0,1\right]  \right)  $ with bounded derivative so that $\varphi_{j}=1$ on
$B(\beta_{1}(j);\beta_{2}(j))$ and $\varphi_{j}=0$ outside of $B(\beta
_{1}(j);2\beta_{2}(j))$.

\medskip

\noindent By Lemma 1 applied to $T_{k}(y)-f(y)$ on $B_{2s_{k}}\cap Y,$ we may
choose $C^{1}$-smooth maps $\delta_{k}:X\rightarrow\mathbb{R}$ so that on each
$B_{2s_{k}}\cap Y$ we have both
\[
\left\vert T_{k}\left(  y\right)  -f\left(  y\right)  -\delta_{k}\left(
y\right)  \right\vert <2^{-k-2}\varepsilon M_{k}^{-1},
\]
and $\left\Vert \delta_{k}^{\prime}\left(  y\right)  \right\Vert _{X^{\ast}%
}<\varepsilon/8$, where $M_{k}=\sum_{i=1}^{k}\widetilde{M}_{i}$ and
$\widetilde{M}_{i}=\sup_{x\in Y\cap B_{2s_{i}}}\left\Vert \varphi_{i}^{\prime
}\left(  x\right)  \right\Vert _{X^{\ast}}$.

\medskip

\noindent Then we also have, for $y\in B_{2s_{k}}\cap Y$ using our estimate above,

\medskip%

\begin{align*}
\left\Vert T_{k}^{\prime}\left(  y\right)  -H(f^{\prime}\left(  y\right)
)-\delta_{k}^{\prime}\left(  y\right)  \right\Vert _{X^{\ast}}  &
\leq\left\Vert T_{k}^{\prime}\left(  y\right)  -H\left(  f^{\prime}\left(
y\right)  \right)  \right\Vert _{X^{\ast}}+\left\Vert \delta_{k}^{\prime
}\left(  y\right)  \right\Vert _{X^{\ast}}\\
& \\
&  <\varepsilon/8C_{0}+\varepsilon/8\leq\varepsilon/4.
\end{align*}

\medskip

Set $\Delta_{i}\left(  x\right)  =T_{k}\left(  x\right)  -\delta_{k}\left(
x\right)  $ if $\beta(i)=B_{s_{k}}$ is a ball from the subcollection
$\{B_{s_{l}}\}_{l=1}^{\infty}$ covering $Y$, and $\Delta_{i}=F\left(
x_{j}\right)  $ if $\beta(i)=B_{r_{j}}$ belongs to the subcollection
$\{B_{r_{l}}\}_{l=1}^{\infty}$ covering $X\setminus Y$.

\medskip

\noindent Next, we define
\[
h_{i}=\varphi_{i}\prod_{k<i}\left(  1-\varphi_{k}\right)  ,
\]
and
\[
g\left(  x\right)  =\sum_{i}h_{i}\left(  x\right)  \Delta_{i}\left(  x\right)
\]

Note that for each $x$, if $n:=n\left(  x\right)  :=\min\left\{  m:x\in
\beta(m)\right\}  $, then because $1-\varphi_{n}\left(  x\right)  =0$ and
$\beta(n)$ is open, it follows from the definition of the $h_{j}$ that there
is a neighbourhood $N\subset\beta(n)$ of $x$ so that for $z\in N$, $g\left(
z\right)  =\sum_{j\leq n}h_{j}\left(  z\right)  \Delta_{j}\left(  z\right)  $,
and $\sum_{j}h_{j}(z)=\sum_{j\leq n}h_{j}(z)$. Also, by a straightforward
calculation, again using the fact that $\varphi_{n}=1$ on $\beta(n)$, we have
that $\sum_{j}h_{j}\left(  z\right)  =1$ for $z\in\beta(n)$, and so for all
$z\in X.$

Now, fix any $x_{0}\in X,$ and let $n_{0}=n\left(  x_{0}\right)  $ and a
neighborhood $N_{0}$ of $x_{0}$ be as above. For each $j\leq n_{0}$ define the
functions $V_{j}:N_{0}\rightarrow\mathbb{R}$ and $W_{j}:N_{0}\rightarrow
\mathbb{R}$ by%
\[
V_{j}(x)=%
\genfrac{\{}{.}{0pt}{}{0\text{ if }\beta_{1}(j)\notin Y}{\left\vert
T_{k}\left(  x\right)  -F\left(  x\right)  -\delta_{k}\left(  x\right)
\right\vert \text{ if }\beta_{1}(j)=y_{k}}%
\]
and%
\[
W_{j}(x)=%
\genfrac{\{}{.}{0pt}{}{0\text{ if }\beta_{1}(j)\in Y}{|F(x_{i})-F(x)|\text{ if
}\beta_{1}(j)=x_{i}}%
\]

\medskip

\noindent Then for any $x\in N_{0}$ we have that
\begin{align*}
\left\vert g\left(  x\right)  -F\left(  x\right)  \right\vert  &  \leq
\sum_{j\leq n_{0}}h_{j}\left(  x\right)  \max\{V_{j}(x),W_{j}(x)\}\\
& \\
&  \leq\sum_{j\leq n_{0}}h_{j}\left(  x\right)  \frac{\varepsilon}%
{2}<\varepsilon.
\end{align*}

\noindent Now define the function $\alpha$ so that when $\beta_{1}(j)\in Y,$
$\beta_{1}(j)=y_{\alpha(j)}$. Recall that $B_{2r_{j}}\cap Y=\emptyset$ for
every $j,$ and that $\varphi_{j}=0$ off of $B_{2r_{j}}.$ Hence, if $y\in Y,$
then in the sum $g\left(  y\right)  $ only those indices $j$ such that
$\beta_{1}(j)\in Y$ are non-zero.

\medskip

\noindent Recall also that $\sum_{j}h_{j}\left(  x\right)  =1$ for all $x\in
X$ (and hence $\sum_{j}h_{j}^{\prime}\left(  x\right)  =0$)$,$ and so for
$y\in Y$ we have, $H\left(  f^{\prime}\left(  y\right)  \right)  =\sum
_{j}h_{j}^{\prime}\left(  y\right)  f\left(  y\right)  +\sum_{j}h_{j}\left(
y\right)  H\left(  f^{\prime}\left(  y\right)  \right)  .$ And, $g^{\prime
}\left(  y\right)  =\sum_{j}h_{j}^{\prime}\left(  y\right)  \left(
T_{\alpha(j)}\left(  y\right)  -\delta_{\alpha(j)}\left(  y\right)  \right)
+\sum_{j}h_{j}\left(  y\right)  \left(  T_{\alpha(j)}^{\prime}\left(
y\right)  -\delta_{\alpha(j)}^{\prime}\left(  y\right)  \right)  .$ Finally, a
straightforward calculation shows that $\left\Vert h_{j}^{\prime
}(x)\right\Vert \leq M_{\alpha(j)}$ for $\beta_{1}(j)\in Y$. With these
observations in mind, we have,

\medskip%

\begin{align*}
&  \left\Vert g^{\prime}\left(  y\right)  -H\left(  f^{\prime}\left(
y\right)  \right)  \right\Vert _{X^{\ast}}\\
& \\
&  \leq\sum_{_{\substack{j\leq n\left(  y\right)  \\\beta_{1}(j)\in Y}%
}}\left(  \left\Vert h_{j}^{\prime}\left(  y\right)  \right\Vert \left\vert
T_{\alpha(j)}\left(  y\right)  -f\left(  y\right)  -\delta_{\alpha(j)}\left(
y\right)  \right\vert \right. \\
&  \left.  +\ h_{j}\left(  y\right)  \left\Vert T_{\alpha(j)}^{\prime}\text{
}(y)-H(f^{\prime}\left(  y\right)  )-\delta_{\alpha(j)}^{\prime}\left(
y\right)  \right\Vert _{X^{\ast}}\right) \\
& \\
&  <\sum_{_{\substack{j\leq n\left(  y\right)  \\\beta_{1}(j)\in Y}%
}}\left\Vert h_{j}^{\prime}\left(  y\right)  \right\Vert \left(
2^{-\alpha(j)-2}\varepsilon M_{\alpha(j)}^{-1}\right)  +\sum
_{_{\substack{j\leq n\left(  y\right)  \\\beta_{1}(j)\in Y}}}h_{j}\left(
y\right)  \frac{\varepsilon}{4}\\
& \\
&  <\sum_{_{\substack{j\leq n\left(  y\right)  \\\beta_{1}(j)\in Y}}}%
M_{\alpha(j)}\left(  2^{-\alpha(j)-2}\varepsilon M_{\alpha(j)}^{-1}\right)
+\frac{\varepsilon}{4}<\varepsilon.
\end{align*}

\medskip

\noindent As $H\left(  f^{\prime}\left(  y\right)  \right)  \mid_{Y}%
=f^{\prime}\left(  y\right)  ,$ we also have the estimate $\left\Vert
g^{\prime}\left(  y\right)  -f^{\prime}\left(  y\right)  \right\Vert
_{Y^{\ast}}<\varepsilon.$

\medskip

\noindent Let us now consider the case when $f$ is $C^{1}$ and Lipschitz on
$Y$ and $F$ is any Lipschitz extension of $f$ to $X$ with $\text{Lip}%
(F)=\text{Lip}(f)$. In this case we have to modify the definition of the
functions $\Delta_{i}$ as follows.

\medskip

\noindent We let $\Delta_{i}\left(  x\right)  =T_{k}\left(  x\right)
-\delta_{k}\left(  x\right)  $ if $\beta(i)=B_{s_{k}}$ is a ball from the
subcollection $\{B_{s_{l}}\}_{l=1}^{\infty}$ covering $Y$ where $\delta_{k}$
is chosen (by using Lemma \ref{Main}) so that $\left\vert T_{k}\left(
y\right)  -f\left(  y\right)  -\delta_{k}\left(  y\right)  \right\vert
<2^{-i-2}\varepsilon M_{k}^{-1}$, $\textrm{Lip}(\delta_k)\leq C_0 \textrm{Lip}(T_k-F)\leq
2C_0\textrm{Lip}(F)$, and $\left\Vert \delta_{k}^{\prime}\left(
y\right)  \right\Vert _{X^{\ast}}<\varepsilon/8$, where now the
$M_{i}$ are
defined by $M_{i}=\sum_{j=1}^{i}\widetilde{M}_{j}$ and $\widetilde{M}_{j}%
=\sup_{x\in B(\beta_{1}(j);2\beta_{2}(j))}\left\Vert
\varphi_{j}^{\prime }\left(  x\right)  \right\Vert _{X^{\ast}}$.
We also let $\Delta _{i}(x)=F_{\ell}(x)$ if
$\beta(i)=B_{r_{\ell}}$ belongs to the subcollection
$\{B_{r_{j}}\}_{j=1}^{\infty}$ covering $X\setminus Y$, where the
function $F_{\ell}$ is again chosen by using Lemma \ref{Main} so
that $\left\vert F_{\ell}(x)-F(x)\right\vert$ $
<2^{-i-2}\varepsilon M_{l}^{-1}$ on $B_{2r_{l}}$ and
$\text{Lip}(F_{\ell})\leq C_{0}\text{Lip}(F)=C_{0}\text{Lip}(f)$.
Note that, with these choices, we have
\[
|\Delta_{i}(x)-F(x)|<2^{-i-2}\varepsilon M_{i}^{-1}\,\,\text{ on }%
\,\,B(\beta_{1}(i),\beta_{2}(i)),\text{ and }%
\]%
\[
\text{Lip}(\Delta_{i})\leq C_{1}\,\text{Lip}(f),
\]
where $C_{1}:=M+3C_0$ (with $M$ as in Lemma \ref{existence of
extension operators}(2)) is a constant depending only on $X$.

\medskip

\noindent Now define the $C^{1}$ function $g:X\rightarrow\mathbb{R}$ by
\[
g(x)=\sum_{i}\Delta_{i}(x)h_{i}(x).
\]
As above, one can check that $\left\vert g(x)-F(x)\right\vert <\varepsilon$
for all $x\in X$, and also $\left\Vert g^{\prime}(y)-H(f^{\prime
}(y))\right\Vert _{X^{\ast}}<\varepsilon$ for all $y\in Y$; that is $g$
satisfies properties $(1)$ and $(2)$ of the statement. Let us see that $g$
satisfies $(3)$ as well. Noting that $\text{Lip}(h_{j})\leq M_{j}$, we can
estimate, for every $x,z\in X$,
\begin{align*}
&  g(x)-g(z)\\
&  =\sum_{j}\Delta_{j}(x)h_{j}(x)-\sum_{j}\Delta_{j}(z)h_{j}(z)\\
&  =\sum_{j}(\Delta_{j}(x)-F(x))(h_{j}(x)-h_{j}(z))+\sum(\Delta_{j}%
(x)-\Delta_{j}(z))h_{j}(z)\\
&  \leq\sum_{j}\frac{\varepsilon}{2^{j+2}M_{j}}\text{Lip}(h_{j})\Vert
x-z\Vert+\sum_{j}\text{Lip}(\Delta_{j})\Vert x-z\Vert h_{j}\left(  z\right) \\
&  \leq\left(  \frac{\varepsilon}{4}+C_{1}\text{Lip}(f)\right)  \Vert
x-z\Vert\leq C\text{Lip}(f)\Vert x-z\Vert,
\end{align*}
provided that $\varepsilon>0$ is chosen small enough (recall that we are
assuming $\text{Lip}(f)>0$), and where $C=2C_{1}>1$, a constant only depending
on $X$. This shows that $\text{Lip}(g)\leq C\,\text{Lip}(f)$. \thinspace
\thinspace\ $\blacksquare$

\medskip

\begin{theorem}
\label{extension theorem} Let $X$ be a separable Banach space which admits a
$C^{1}$-smooth norm. Let $Y\subset X$ be a closed subspace, and
$f:Y\rightarrow\mathbb{R}$ a $C^{1}$-smooth and Lipschitz function. Then there
is a $C^{1}$ and Lipschitz extension $g:X\to\mathbb{R}$ of $f$ such that
$\text{Lip}(g)\leq C\text{Lip}(f)$, where $C$ is a constant depending only on
$X$.
\end{theorem}

\noindent\textbf{Proof.\ \ } First note that if $h$ is a bounded, Lipschitz
function defined on $Y$, there always exists a bounded, Lipschitz extension of
$h$ to $X$, with the same Lipschitz constant, and bounded by the same constant
(defined for instance by $x\mapsto\max\{ -\|h\|_{\infty}, \min\{
\|h\|_{\infty}, \, \inf_{y\in Y}\{ h(y)+\text{Lip}(h)\|x-y\|\}\, \}\, \}$).
For the purposes of the proof, we denote such an extension by $\overline{h}$.

\medskip

\noindent We are going to define our function $g$ by means of a series
constructed by induction. By Theorem 1 there exists a $C^{1}$ function
$g_{1}:X\rightarrow\mathbb{R}$ such that

\begin{itemize}
\item $|\overline{f-g_{1}}|<2^{-1}\varepsilon$ on $X$,

\item $\left\Vert f^{\prime}\left(  y\right)  -g_{1}^{\prime}\left(  y\right)
\right\Vert _{Y^{\ast}}<2^{-1}\varepsilon/C$ for $y\in Y$ (note in particular
that this implies $\text{Lip}(f-g_{|_{Y}})\leq2^{-1}\varepsilon/C$), and

\item $\text{Lip}\left(  g_{1}\right)  \leq C\,\text{Lip}(f)$
\end{itemize}

\medskip

\noindent Now, for $n\geq2$, suppose that we have chosen $g_{1},...,g_{n}%
,\;$real-valued and $C^{1}$-smooth on $X$ such that for all $x\in X$ and $y\in
Y,$

\medskip%

\[
\left\vert \overline{\left(  f-\sum_{i=1}^{n}g_{i}\right)  }\left(  x\right)
\right\vert <2^{-n}\varepsilon,
\]

\medskip%

\[
\left\Vert f^{\prime}\left(  y\right)  -\sum_{i=1}^{n}g_{i}^{\prime}\left(
y\right)  \right\Vert _{Y^{\ast}}<2^{-n}\varepsilon/C,
\]

\medskip

\noindent and%

\[
\text{Lip}(g_{n})\leq C\ \text{Lip}\left(  f-(\sum_{j=1}^{n-1}{g_{j}}%
)_{|Y}\right)  .
\]
It is clear that an application of Theorem 1 to the function $f-(g_{1}%
)_{|_{Y}}$ provides us with a function $g_{2}$ which, together with $g_{1}$,
makes the above properties true for $n=2$. Hence we can proceed to the general
step of our inductive construction.

\medskip

\noindent Consider the function $l=f-\sum_{i=1}^{n}g_{i},$ which is $C^{1}%
$-smooth on $Y.$

\medskip

\noindent By Theorem 1, we can find a $C^{1}$-smooth map $g_{n+1}$ on $X$ such
that we have,

\medskip%

\begin{equation}
\left\vert \overline{l}\left(  x\right)  -g_{n+1}\left(  x\right)  \right\vert
<2^{-n-1}\varepsilon\,\,\text{ on\ }X,
\end{equation}
\noindent and for $y\in Y,$
\begin{equation}
\left\Vert l^{\prime}\left(  y\right)  -g_{n+1}^{\prime}\left(  y\right)
\right\Vert _{Y^{\ast}}<2^{-n-1}\varepsilon/C,
\end{equation}

\noindent and also,%

\begin{equation}
\text{Lip}\left(  g_{n+1}\right)  \leq C\,\text{Lip}\left(  l\right)  .
\end{equation}

\medskip

\noindent From $\left(  2.5\right)  ,$ we have in particular, $\left\vert
f\left(  y\right)  -\sum_{i=1}^{n+1}g_{i}\left(  y\right)  \right\vert
=\left\vert l\left(  y\right)  -g_{n+1}\left(  y\right)  \right\vert
<2^{-n-1}\varepsilon\ $on\ $Y,$ and so $\left\vert \overline{\left(
f-\sum_{i=1}^{n+1}g_{i}\right)  }\left(  x\right)  \right\vert <2^{-n-1}%
\varepsilon$ on $X.$ This together with $\left(  2.6\right)  $ and $\left(
2.7\right)  $ completes the inductive step.

\medskip

\noindent Now, from $\left(  2.6\right)  $ we have $\text{Lip}\left(
f-(\sum_{i=1}^{n+1}g_{i})_{|Y}\right)  \leq2^{-n}\varepsilon/C$, and so from
$\left(  2.7\right)  $ we obtain,%

\[
\Vert g_{n+1}^{\prime}(x)\Vert\leq\text{Lip}(g_{n+1})\leq C\,\text{Lip}\left(
f-(\sum_{j=1}^{n}g_{j})_{|Y}\right)  \leq C2^{-n}\varepsilon/C=2^{-n}%
\varepsilon.
\]
Hence the series $\sum_{j}g_{j}^{\prime}(x)$ is absolutely and uniformly
convergent on $X$. Similarly, we have the estimate $|g_{n+1}(x)|\leq
2^{-n+1}\varepsilon$. Therefore the series
\[
g(x)=\sum_{n=1}^{\infty}g_{n}(x)
\]
defines a $C^{1}$ function on $X$, which coincides with $f$ on $Y$ because of
the first inequality in the inductive assumptions. Finally, we have
\[
\text{Lip}(g)\leq\text{Lip}(g_{1})+\sum_{n=2}^{\infty}\text{Lip}(g_{n})\leq
C\ \text{Lip}(f)+\sum_{n=2}^{\infty}2^{-\left(  n-1\right)  }\varepsilon
\leq2C\ \text{Lip}(f),
\]
provided that $\text{Lip}(f)>0$ (which we can always assume) and $\varepsilon$
is small enough. $\blacksquare$

\medskip

\begin{remark}
\emph{With some more work in the proofs of the preceding theorems one could
show that the constant $C$ can be taken to be any number $C>M$, where $M$ is
as in Lemma \ref{existence of extension operators}(2). Unfortunately the proof
of the Bartle-Graves extension theorem does not give us any useful estimation
about the size of $M$, and in general $M$ is going to be quite large, so we
cannot hope that any refinement of the above proofs will yield a statement of
Theorem \ref{extension theorem} in which $C$ can be chosen to be any number
bigger than $1$.}
\end{remark}

\medskip

\begin{theorem}
Let $X$ be a separable Banach space which admits a $C^{1}$-smooth norm. Let
$Y\subset X$ be a closed subspace, and $f:Y\rightarrow\mathbb{R}$ a $C^{1}%
$-smooth function. Then there is a $C^{1}$ extension of $f$ to $X$.
\end{theorem}

\noindent\textbf{Proof.\ \ } Since $f$ is $C^{1}$ on $Y$, there exists
$\{B_{j}\}:=\{B(y_{j};r_{j})\}_{j=1}^{\infty}$, a countable covering of $Y$ by
open balls in $X$ such that $f$ is Lipschitz on $B_{j}\cap Y$ for each
$j\in\mathbb{N}$. Let $U=\bigcup_{j=1}^{\infty}B_{j}$, and $W=X\setminus Y$.

Consider the mapping $h:X\rightarrow X$ defined by
\[
h(x)=\frac{1}{1+\Vert x\Vert}\,x.
\]
It is easily checked that $h$ is a $C^{1}$ diffeomorphism from $X$ onto its
open unit ball $\mathbf{int}B_{X}$ (with inverse $h^{-1}(y)=\left(  1/(1-\Vert
y\Vert)\right)  y$), that $h$ has a bounded derivative, and that $h$ preserves
lines and in particular leaves the subspace $Y$ invariant. By composing $h$
with suitable dilations and translations we get $C^{1}$ diffeomorphisms
$h_{j}:X\rightarrow B_{j}$ such that $h_{j}$ is Lipschitz for each $j$. And,
by composing the restrictions to $Y$ of these $h_{j}$ with our function $f$,
we get $C^{1}$ and Lipschitz functions $f_{j}:=f\circ(h_{j})_{|_{Y}%
}:Y\rightarrow\mathbb{R}$. According to the preceding result there exist
$C^{1}$ (and Lipschitz) extensions $G_{j}:X\rightarrow\mathbb{R}$ of $f_{j}$.
Then the composition
\[
g_{j}=G_{j}\circ h_{j}^{-1}%
\]
defines a $C^{1}$ extension of $f_{|_{B_{j}\cap Y}}$ to $B_{j}$. Put
$g_{0}\equiv1.$

Now let $\{\varphi_{0}\}\cup\{\varphi_{j}\}_{j=1}^{\infty}$ be a $C^{1}$
partition of unity subordinated to the open covering $\{W\}\cup\{B_{j}%
\}_{j=1}^{\infty}$ of $X$ (such partitions of unity always exist for separable
spaces with $C^{1}$ norms, see \cite[Theorem VIII.3.2, page 351]{DGZ}).
Define
\[
g(x)=\sum_{j=0}^{\infty}\varphi_{j}(x)g_{j}(x).
\]
Then it is clear that $g$ is a $C^{1}$ extension of $f$ to $X$. $\blacksquare$

\bigskip

\begin{corollary}
Let $M$ be a separable Banach manifold modelled on a Banach space $X$ which
admits a $C^{1}$ norm, and let $N$ be a closed $C^{1}$ submanifold of $M$.
Then every $C^{1}$ function $f:N\to\mathbb{R}$ has a $C^{1}$ extension to $M$.
\end{corollary}

\noindent\textbf{Proof.\ \ } Let $\{V_{j}\}_{j=1}^{\infty}$ be a covering of
$N$ by open sets in $M$ so that there are $C^{1}$ diffeomorphisms $\psi
_{j}:V_{j}\to X$ such that $\psi_{j} (N\cap V_{j})=Y$, where $Y$ is a closed
subspace of $X$.

The functions $f\circ\psi_{j}^{-1}:Y\to\mathbb{R}$ are $C^{1}$ and (by the
preceding theorem) there are $C^{1}$ extensions $G_{j}:X\to\mathbb{R}$, which
in turn give, by composition, $C^{1}$ extensions $g_{j}:=G_{j}\circ\psi_{j}$
of $f_{|_{V_{j}\cap N}}$ to $V_{j}$.

Then, if $\{\theta\}\cup\{\theta_{j}\}_{j=1}^{\infty}$ is a $C^{1}$ partition
of unity subordinated to the open covering $\{M\setminus N\}\cup
\{V_{j}\}_{j=1}^{\infty}$ of $M$ (note that a separable Banach manifold
modelled on a Banach space $X$ admits $C^{1}$ partitions of unity if and only
if $X$ does), the function
\[
g(x)=\sum_{j}\theta_{j}(x)g_{j}(x)
\]
is a $C^{1}$ extension of $f$ to $M$. $\blacksquare$

\medskip

If $Y$ is not required to be a closed subspace of $X$ but is
merely closed, results similar to Theorem 1 and Theorem 2 can be
obtained. However, the differentiability requirements on $f$ must
be strengthened. The proofs, which we omit, closely parallel those
for Theorem 1 and Theorem 2, where the essential difference is
that $f^{\prime}(y)$ is extended to directions off of $Y$, not by
the Bartle-Graves generated $H(f^{\prime}(y))$, but by explicit
hypothesis. Of course one must also verify that Lemma 1 still
holds in the case when $Y$ is a closed subset of $X$. It is easy
to establish such a version of Lemma 1 by replacing the function
$\|q(x)\|$ in its proof with a Lipschitz $C^1$ approximation of
the distance function to $Y$ (which in turn can be constructed
with the help of a sup-partition of unity provided by Claim 1, see
also \cite{F, HJ}).

\medskip

\begin{theorem}
\label{Key}Let $X$ be a separable Banach space which admits a $C^{1}$-smooth
norm, $Y\subset X$ a closed subset, and $U\supset Y$ a neighbourhood of $Y.$
Let $\varepsilon>0,$ and $f:U\rightarrow\mathbb{R}$ be a map which is $C^{1}%
$-smooth on $Y$ as a function on $X.$ Then there exists a $C^{1}$-smooth map
$g:X\rightarrow\mathbb{R}$ such that,

\begin{enumerate}
\item $\left\vert f\left(  y\right)  -g\left(  y\right)  \right\vert
<\varepsilon$ on $Y,$

\smallskip

\item $\left\Vert f^{\prime}\left(  y\right)  -g^{\prime}\left(  y\right)
\right\Vert _{X^{\ast}}<\varepsilon$ on $Y.$
\end{enumerate}
\end{theorem}

\medskip

\begin{theorem}
\label{Mainresult} Let $X$ be a separable Banach space which admits a $C^{1}%
$-smooth norm. Let $Y\subset X$ be a closed subset, and $U$ an open set
containing $Y.$ Let $f:U\rightarrow\mathbb{R}$ be a map $C^{1}$-smooth on $Y$
as a function on $X.$ Then there is a $C^{1}$ extension of $f\mid_{Y}$ to $X$.
\end{theorem}

\noindent We have the following easy corollary.

\begin{corollary}
Let $X$ be a separable Banach space which admits a $C^{1}$-smooth norm. Let
$U\subset X$ be open and $f:U\rightarrow\mathbb{R}$ a $C^{1}$-smooth function.
Then for any open set $V\subset U$ with $\overline{V}\subset U,$ there is a
$C^{1}$ extension of $f\mid_{\overline{V}}$ to $X$.
\end{corollary}


\medskip

\noindent\textbf{Acknowledgement.} The authors wish to thank
Richard Aron, Gilles Godefroy, and Sophie Grivaux for many helpful
discussions. We also thank P. Hajek and M. Johannis for drawing
our attention to their recent work [HJ], which led to a
significant improvement of the first version of our paper, and to
M. Jim{\'e}nez-Sevilla and L. S{\'a}nchez-Gonz{\'a}lez for pointing out a gap
in the proof of Theorem 1 which this new version fixes.

\end{document}